\documentclass{amsart}

\usepackage{a4wide,latexsym,amssymb,amsthm,amsmath}

\theoremstyle{plain}

\newtheorem{theorem}{Theorem}
\newtheorem{lemma}{Lemma}
\newtheorem{proposition}{Proposition}

\theoremstyle{remark}

\newtheorem{remark}{Remark}
\newtheorem{example}{Example}

\def\nt{\widetilde{\nabla}}
\def\M{\widetilde{M}}

\def\R{\mathbb{R}}
\def\C{\mathbb{C}}
\def\E{\mathbb{E}}
\def\H{\mathbb{H}}
\def\S{\mathbb{S}}
\def\D{\mathbb{D}}
\def\L{\mathbb{L}}

\def\e{\eqref}

\def\<{ \left < }
\def\>{ \right > }

\newcommand{\Nil}{\mathrm{Nil}_3}
\newcommand{\Sol}{\mathrm{Sol}_3}

\begin{document}

\title[Higher order parallel surfaces in Bianchi-Cartan-Vranceanu spaces]
{Higher order parallel surfaces in Bianchi-Cartan-Vranceanu spaces}

\author[J. Van der Veken]{Joeri Van der Veken}
\address{Katholieke Universiteit Leuven\\ Departement
Wiskunde\\ Celestijnenlaan 200 B\\ B-3001 Leuven\\ Belgium}
\email{joeri.vanderveken@wis.kuleuven.be}
\thanks{The author is a postdoctoral researcher supported by the Research Foundation-Flanders (FWO)}

\begin{abstract}
We give a full classification of higher order parallel surfaces in
three-dimen\-sional homogeneous spaces with four-dimensional
isometry group, i.e. in the so-called Bianchi-Cartan-Vranceanu
family. This gives a positive answer to a conjecture formulated in
\cite{3}. As a partial result, we prove that totally umbilical
surfaces only exist if the ambient Bianchi-Cartan-Vranceanu space is
a Riemannian product of a surface of constant Gaussian curvature and
the real line, and we give a local parametrization of all totally
umbilical surfaces.
\end{abstract}



\maketitle

\section{Introduction}

A Riemannian manifold $(M,g)$ is said to be homogeneous if for every
two points $p$ and $q$ in $M$, there exists an isometry of $M$,
mapping $p$ into $q$. The classification of simply connected
3-dimensional homogeneous spaces is well-known. The dimension of the
isometry group must equal 6, 4 or 3. If the isometry group is of
dimension 6, $M$ is a complete real space form, i.e. Euclidean space
$\E^3$, a sphere $\S^3(\kappa)$, or a hyperbolic space
$\H^3(\kappa)$. If the dimension of the isometry group is 4, $M$ is
isometric to $\mathrm{SU}(2)$, the special unitary group, to
$[\mathrm{SL}(2,\R)]^{\sim}$, the universal covering of the real
special linear group, to $\Nil$, the Heisenberg group, all with a
certain left-invariant metric, or to a Riemannian product
$\S^2(\kappa)\times\R$ or $\H^2(\kappa)\times\R$. Finally, if the
dimension of the isometry group is 3, $M$ is isometric to a general
simply connected Lie group with left-invariant metric. As will
become clear in the next section, Bianchi-Cartan-Vranceanu spaces
are in fact the spaces with 4-dimensional isometry group mentioned
above, together with $\E^3$ and $\S^3(\kappa)$.\medskip

The classification above contains the eight ``model geometries''
appearing in the famous conjecture of Thurston on the classification
of 3-manifolds, namely $\E^3$, $\S^3$, $\H^3$, $\S^2\times\R$,
$\H^2\times\R$, $[\mathrm{SL}(2,\R)]^{\sim}$, $\Nil$ and $\Sol$. See
for example \cite{17a}. In theoretical cosmology, the metrics on
these spaces are known as Bianchi-Kantowski-Sachs type metrics, used
to construct spatially homogeneous spacetimes, see for example
\cite{13a}.\medskip

Immersions of curves and surfaces in 3-dimensional real space forms
are extensively studied and it is now very natural to allow the
other 3-dimensional homogeneous manifolds as ambient spaces. Initial
work in this direction can be found in \cite{12} and
\cite{7}.\medskip

An important class of surfaces to study are parallel surfaces. These
immersions have a parallel second fundamental form and hence their
extrinsic invariants ``are the same'' at every point. Parallel
submanifolds in real space forms are classified in \cite{1}. In
\cite{8}, \cite{9}, \cite{10} and \cite{15}, the notion of higher
order parallelism is introduced and a classification for
hypersurfaces in real space forms is obtained. In \cite{3}, a
classification of parallel surfaces in 3-dimensional homogeneous
spaces with 4-dimensional isometry group is given, whereas the
classification of higher order parallel surfaces is formulated as a
conjecture. In this article we will prove this conjecture (Theorem
\ref{theo7}). For an overview of the theory of parallel and higher
order parallel submanifolds we refer to \cite{16}.\medskip

Another important class of surfaces are totally umbilical ones. From
an extrinsic viewpoint, these surfaces are curved equally in every
direction. We will give a full local classification of totally
umbilical surfaces in 3-dimensional homogeneous spaces with
4-dimensional isometry group (Theorems \ref{theo}, \ref{T:3.3.2} and
\ref{T:3.3.3}). Although in a real space form a totally umbilical
surface is automatically parallel, this will no longer be the case
in the
spaces under consideration.\\

\section{Examples of three-dimensional homogeneous spaces}

\subsection{The Heisenberg group $\Nil$ with left-invariant metric}
The Heisenberg group $\Nil$ is  a Lie group which is diffeomorphic
to $\R^3$ and the group operation is defined by
$$(x,y,z)\ast(\overline{x},\overline{y},\overline{z})=\left(x+\overline{x},\ y+\overline{y},\ z+\overline{z}+\frac{x\overline{y}}{2}-\frac{\overline{x}y}{2}\right).$$
Remark that the mapping
$$\Nil\rightarrow\left\{\left.\left(\begin{array}{ccc}1&a&b\\0&1&c\\0&0&1\end{array}\right)\ \right|\ a,b,c\in\R\right\} : (x,y,z)\mapsto\left(\begin{array}{ccc}1&x&z+\frac{xy}{2}\\0&1&y\\0&0&1\end{array}\right)$$
is an isomorphism between $\Nil$ and a subgroup of
$\mathrm{GL}(3,\R)$. For every non-zero real number $\tau$ the
following metric on $\Nil$ is left-invariant:
$$ds^2=dx^2+dy^2+4\tau^2\left(dz+\frac{y\,dx-x\,dy}{2}\right)^2.$$
After the change of coordinates $(x,y,2\tau z)\mapsto(x,y,z)$, this
metric is expressed as
\begin{equation} \label{2.1}
ds^2=dx^2+dy^2+\left(dz+\tau(y\,dx-x\,dy)\right)^2.
\end{equation}
\medskip

\subsection{The projective special linear group $\mathrm{PSL}(2,\R)$ with left-invariant metric}
Consider the following subgroup of $\mathrm{GL}(2,\R)$:
$$\mathrm{SL}(2,\R)=\left\{\left.\left(\begin{array}{cc}a&b\\c&d\end{array}\right)\ \right|\ ad-bc=1\right\}.$$
First note that this group is isomorphic to the
following subgroup of $\mathrm{GL}(2,\C)$:
$$G=\left\{\left.\left(\begin{array}{cc}\alpha&\beta\\\overline{\beta}&\overline{\alpha}\end{array}\right)\ \right|\ |\alpha|^2-|\beta|^2=1\right\},$$
via the isomorphism
$$\mathrm{SL}(2,\R)\rightarrow
G:\left(\begin{array}{cc}a&b\\c&d\end{array}\right)\mapsto\frac{1}{2}\left(\begin{array}{cc}i&1\\1&i\end{array}\right)\left(\begin{array}{cc}a&b\\c&d\end{array}\right)\left(\begin{array}{cc}-i&1\\1&-i\end{array}\right).$$
Now consider the Poincar\'e disc-model for the hyperbolic plane
$\H^2(\kappa)$ of constant Gaussian curvature $\kappa<0$:
\begin{eqnarray} \label{2.2}
\H^2(\kappa) &\cong& \left(\left\{(x,y)\in\R^2\ \left|\ x^2+y^2<-\frac{4}{\kappa}\right.\right\},\ \frac{dx^2+dy^2}{(1+\frac{\kappa}{4}(x^2+y^2))^2}\right)\\
             &\cong& \left(\left\{z\in\C\ \left|\ |z|^2<-\frac{4}{\kappa}\right.\right\},\ \frac{dz\,d\overline{z}}{(1+\frac{\kappa}{4}|z|^2)^2}\right)\nonumber
\end{eqnarray}
and define
$$F_{\left(\alpha\ \beta\atop\overline{\beta}\ \overline{\alpha}\right)}(z)=\frac{2}{\sqrt{-\kappa}}\frac{\alpha\sqrt{-\kappa}z + 2\beta}{\overline{\beta}\sqrt{-\kappa}z + 2\overline{\alpha}}.$$
Note that for $\kappa=-4$, this M\"obius transformation simplifies
to $z\mapsto\frac{\alpha z+\beta}{\overline{\beta}
z+\overline{\alpha}}$. The mapping
$$G\times\H^2(\kappa)\rightarrow\H^2(\kappa):(A,z)\mapsto F_A(z)$$
is a transitive, isometric action with stabilizers isomorphic to the
circle group $\mathrm{SU}(1)=\{\alpha\in\C\ |\ |\alpha|=1\}.$ This
action induces the following transitive action on the unitary
tangent bundle $\mathrm{U}\H^2(\kappa)$:
\begin{equation} \label{2.3}
G\times\mathrm{U}\H^2(\kappa)\rightarrow\mathrm{U}\H^2(\kappa):(A,(z,v))\mapsto(F_A(z),(F_A)_{\ast}v),
\end{equation}
with stabilizers of order two. Hence we can identify
$\mathrm{U}\H^2(\kappa)$ with
$$\mathrm{PSL}(2,\R)=\frac{\mathrm{SL}(2,\R)}{\left\{\left(1\ 0\atop 0\ 1\right),-\left(1\ 0\atop 0\ 1\right)\right\}}.$$

Let us now define a metric on $\mathrm{U}\H^2(\kappa)$. If
$\gamma:I\subseteq\R\rightarrow\mathrm{U}\H^2(\kappa):t\mapsto
(z(t),v(t))$ is a curve, with $z(t)$ a curve in $\H^2(\kappa)$ and
for every $t\in I$, $v(t)\in T_{z(t)}\H^2(\kappa)$ and $\|v(t)\|=1$,
we put
\begin{equation} \label{2.4}
\|\gamma'(t_0)\|^2 = \|z'(t_0)\|^2 +
\left(\frac{2\tau}{\kappa}\right)^2\|(\nabla_{z'}v)_{z(t_0)}\|^2,\quad
\tau\in\R\setminus\{0\},
\end{equation}
where $\nabla_{z'}v$ is the covariant derivative of the vector field
$v$ along the curve $z(t)$. For $\tau=\pm\frac{\kappa}{2}$, this
metric is induced from the standard metric on the tangent bundle. By
varying the parameter $\tau$, we distort the length of the fibres.
It is clear that the action (\ref{2.3}) is now isometric and hence
the induced metric on $\mathrm{PSL}(2,\R)$ via the identification is
left-invariant. The metric (\ref{2.4}) can be explicitly computed,
analogous as in \cite{7}, in the coordinate system
\begin{multline*}
\D^2\left(\frac{2}{\sqrt{-\kappa}}\right)\times\S^1(1)\rightarrow\mathrm{U}\H^2(\kappa):\\
((x,y),\theta)\mapsto\left((x,y),\
\left(1+\frac{\kappa}{4}(x^2+y^2)\right)\left(\cos\left(\frac{\kappa}{2\tau}\theta\right)\frac{\partial}{\partial
x}+\sin\left(\frac{\kappa}{2\tau}\theta\right)\frac{\partial}{\partial
y}\right)\right),
\end{multline*}
where $\D^2\left(\frac{2}{\sqrt{-\kappa}}\right)$ is the disc of
radius $\frac{2}{\sqrt{-\kappa}}$, yielding the following result:
\begin{equation} \label{2.5}
ds^2=\frac{dx^2+dy^2}{(1+\frac{\kappa}{4}(x^2+y^2))^2}+\left(d\theta+\tau\frac{y\,dx-x\,dy}{1+\frac{\kappa}{4}(x^2+y^2)}\right)^2,\quad
\kappa<0.
\end{equation}
\medskip

\subsection{The special orthogonal group $\mathrm{SO}(3)$ with left-invariant metric}
Consider the following subgroup of $\mathrm{GL}(2,\C)$:
$$\mathrm{SU}(2)=\left\{\left.\left(\begin{array}{cc}\alpha&\beta\\-\overline{\beta}&\overline{\alpha}\end{array}\right)\ \right|\ |\alpha|^2+|\beta|^2=1\right\}.$$
Using stereographic projection, we have for an arbitrary
$\kappa>0$:
\begin{equation} \label{2.8}
\S^2(\kappa)\setminus\{\infty\}\cong\left(\R^2,\
\frac{dx^2+dy^2}{(1+\frac{\kappa}{4}(x^2+y^2))^2}\right)\cong\left(\C,\
\frac{dz\,d\overline{z}}{(1+\frac{\kappa}{4}|z|^2)^2}\right).
\end{equation}
The analogy with the previous case is clear and we could now
proceed in the same way as above, putting
$$F_{\left(\alpha\ \ \beta\atop-\overline{\beta}\ \overline{\alpha}\right)}(z)=\frac{2}{\sqrt{\kappa}}\frac{\alpha\sqrt{\kappa}z + 2\beta}{(-\overline{\beta}\sqrt{\kappa}z + 2\overline{\alpha})},$$
and being careful in calculations involving the symbol $\infty$. In
this way we would find that $\mathrm{U}\S^2(\kappa)$ can be
identified with
$$\mathrm{PSU}(2)=\frac{\mathrm{SU}(2)}{\left\{\left(1\ 0\atop 0\ 1\right),-\left(1\ 0\atop 0\ 1\right)\right\}}.$$
But since $\mathrm{PSU}(2)$ is isomorphic to $\mathrm{SO}(3)$, see
for example \cite{18}, there is an easier way to construct the
desired group action. Looking at $\S^2(\kappa)$ as a hypersphere in
$\E^3$ centered at the origin, we can identify both points of the
surface and tangent vectors to it with elements of $\R^3$ and we
define
$$\mathrm{SO}(3)\times\mathrm{U}\S^2(\kappa)\rightarrow\mathrm{U}\S^2(\kappa):(A,(p,v))\mapsto (Ap,Av).$$
This is a transitive action with trivial stabilizers and a metric on
$\mathrm{U}\S^2(\kappa)$ analogous to (\ref{2.4}) turns it into an
isometric action. This means that the induced metric on
$\mathrm{SO}(3)$ will be left-invariant and in the local coordinates
$$\R^2\times\S^1(1)\rightarrow\mathrm{U}\S^2(\kappa):((x,y),\theta)\mapsto\left((x,y),\ \left(1+\frac{\kappa}{4}(x^2+y^2)\right)\left(\cos\left(\frac{\kappa}{2\tau}\theta\right)\frac{\partial}{\partial x}+\sin\left(\frac{\kappa}{2\tau}\theta\right)\frac{\partial}{\partial y}\right)\right)$$
it is expressed as
\begin{equation} \label{2.9}
ds^2=\frac{dx^2+dy^2}{(1+\frac{\kappa}{4}(x^2+y^2))^2}+\left(d\theta+\tau\frac{y\,dx-x\,dy}{1+\frac{\kappa}{4}(x^2+y^2)}\right)^2,\quad
\kappa>0.
\end{equation}
\medskip

\subsection{The Riemannian product spaces $\H^2(\kappa)\times\R$ and $\S^2(\kappa)\times\R$}Using respectively
the models (\ref{2.2}) and (\ref{2.8}) for $\H^2(\kappa)$ and
$\S^2(\kappa)$, one sees that the Riemannian product metric on these
spaces can be expressed (locally) as
\begin{equation}\label{2.10}
ds^2 = \frac{dx^2+dy^2}{(1+\frac{\kappa}{4}(x^2+y^2))^2}+ dz^2.
\end{equation}
\medskip

\subsection{Bianchi-Cartan-Vranceanu spaces}
Remark that the metrics (\ref{2.1}), (\ref{2.5}), (\ref{2.9}) and
(\ref{2.10}) of the homogeneous spaces above are of the same type.
Cartan classified all 3-dimensional spaces with 4-dimensional
isometry group in \cite{6}. In particular, he proved that they are
all homogeneous and obtained the following two-parameter family of
spaces, which are now known as the {\em Bianchi-Cartan-Vranceanu
spaces} or {\em BCV spaces} for short. For $\kappa,\tau\in\R$, we
define $\M^3(\kappa,\tau)$ as the following open subset of $\R^3$:
$$\left\{(x,y,z)\in\R^3\ \left|\ 1+\frac{\kappa}{4}(x^2+y^2)>0\right.\right\},$$
equipped with the metric
\begin{equation}\label{2.11}
ds^2=\frac{dx^2+dy^2}{(1+\frac{\kappa}{4}(x^2+y^2))^2}+\left(dz+\tau\frac{y\,dx-x\,dy}{1+\frac{\kappa}{4}(x^2+y^2)}\right)^2.
\end{equation}
See also \cite{4}, \cite{5} and \cite{19}. The result of Cartan
shows that the examples above cover in fact all possible
3-dimensional homogeneous spaces with 4-dimensional isometry group.
The BCV family also includes two real space forms, which have
6-dimensional isometry group. The full classification of these
spaces is as follows:
\begin{itemize}
\item if $\kappa=\tau=0$, then $\widetilde{M}^3(\kappa,\tau)\cong \E^3$;
\item if $\kappa=4\tau^2\neq 0$, then $\widetilde{M}^3(\kappa,\tau)\cong \S^3\left(\frac{\kappa}{4}\right)\setminus\{\infty\}$;
\item if $\kappa >0$ and $\tau=0$, then $\widetilde{M}^3(\kappa,\tau)\cong (\S^2(\kappa)\setminus\{\infty\})\times\R$;
\item if $\kappa <0$ and $\tau=0$, then $\widetilde{M}^3(\kappa,\tau)\cong \H^2(\kappa)\times\R$;
\item if $\kappa >0$ and $\tau\neq 0$, then $\widetilde{M}^3(\kappa,\tau)\cong [U(\S^2(\kappa)\setminus\{\infty\})]^{\sim}\cong\mathrm{SU}(2)\setminus\{\infty\}$;
\item if $\kappa <0$ and $\tau\neq 0$, then $\widetilde{M}^3(\kappa,\tau)\cong [U\H^2(\kappa)]^{\sim}\cong[\mathrm{SL}(2,\R)]^{\sim}$;
\item if $\kappa =0$ and $\tau\neq 0$, then $\widetilde{M}^3(\kappa,\tau)\cong \Nil$.
\end{itemize}

To end this section, we discuss the geometry of these spaces. The
following vector fields form an orthonormal frame on
$\M^3(\kappa,\tau)$:
\begin{eqnarray*}
e_1 =
\left(1+\frac{\kappa}{4}(x^2+y^2)\right)\frac{\partial}{\partial
x}-\tau\frac{\partial}{\partial z},\quad e_2 =
\left(1+\frac{\kappa}{4}(x^2+y^2)\right)\frac{\partial}{\partial
y}+\tau\frac{\partial}{\partial z},\quad e_3 =
\frac{\partial}{\partial z}.
\end{eqnarray*}
It is clear that these vector fields satisfy the following
commutation relations:
\begin{equation}\label{Lie-haken}
[e_1,e_2]=-\frac{\kappa}{2}ye_1+\frac{\kappa}{2}xe_2+2\tau e_3,
\qquad [e_2,e_3]=0, \qquad [e_3,e_1]=0.
\end{equation}
The Levi Civita connection of $\M^3(\kappa,\tau)$ can then be
computed using Koszul's formula:
\begin{equation} \label{Levi Civita}
\begin{array}{ccc} \widetilde{\nabla}_{e_1}e_1=\frac{\kappa}{2}ye_2, & \widetilde{\nabla}_{e_1}e_2=-\frac{\kappa}{2}ye_1+\tau e_3, & \widetilde{\nabla}_{e_1}e_3=-\tau e_2, \\
                   \widetilde{\nabla}_{e_2}e_1=-\frac{\kappa}{2}xe_2-\tau e_3, & \widetilde{\nabla}_{e_2}e_2=\frac{\kappa}{2}xe_1, & \widetilde{\nabla}_{e_2}e_3=\tau e_1, \\
                   \widetilde{\nabla}_{e_3}e_1=-\tau e_2, & \widetilde{\nabla}_{e_3}e_2=\tau e_1, & \widetilde{\nabla}_{e_3}e_3=0. \end{array}
\end{equation}
Remark that $\nt_Xe_3=\tau(X\times e_3)$ for every $X\in
T\M^3(\kappa,\tau)$, where the cross product is defined as an
anti-symmetric bilinear operation, satisfying $e_1\times e_2=e_3$,
$e_2\times e_3=e_1$ and $e_3\times e_1=e_2$. The equations in
(\ref{Levi Civita}) yield the following expression for the curvature
tensor of $\M^3(\kappa,\tau)$:
\begin{multline}
\widetilde{R}(X,Y)Z=(\kappa-3\tau^2)(\langle Y,Z\rangle X-\langle X,Z\rangle Y)\\
    -(\kappa-4\tau^2)(\langle Y,e_3 \rangle\langle Z,e_3 \rangle X - \langle X,e_3 \rangle\langle Z,e_3 \rangle Y
    + \langle X,e_3 \rangle\langle Y,Z \rangle e_3 - \langle Y,e_3 \rangle\langle X,Z \rangle e_3)
\end{multline}
for $p\in \M^3(\kappa,\tau)$ and $X,Y,Z\in
T_p\M^3(\kappa,\tau)$.\medskip

Consider the following Riemannian surface with constant Gaussian
curvature $\kappa$:
$$\M^2(\kappa)=\left(\left\{(x,y)\in\R^2\ \left|\ 1+\frac{\kappa}{4}(x^2+y^2)>0\right.\right\}\ ,\ \frac{dx^2+dy^2}{(1+\frac{\kappa}{4}(x^2+y^2))^2}\right).$$
Then the mapping
$$\pi:\M^3(\kappa,\tau)\rightarrow\M^2(\kappa):(x,y,z)\mapsto (x,y)$$
is a Riemannian submersion, referred to as the
\emph{Hopf-fibration}. For $\kappa=4\tau^2\neq 0$, this mapping
coincides with the ``classical'' Hopf-fibration
$\pi:\S^3\left(\frac{\kappa}{4}\right)\rightarrow\S^2(\kappa)$. In
the following, by a \emph{Hopf-cylinder} we mean the inverse image
of a curve in $\M^2(\kappa)$ under $\pi$. By a \emph{leaf} of the
Hopf-fibration, we mean a surface which is everywhere orthogonal to
the fibres. From Frobenius' theorem and (\ref{Lie-haken}), it is
clear that this only exists if $\tau=0$.

\section{Surfaces immersed in BCV spaces}

Let us start with recalling the basic formulas from the theory of
submanifolds. Suppose that $F:M^n\rightarrow\M^{n+k}$ is an
isometric immersion of Riemannian manifolds and denote by $\nabla$
the Levi Civita connection of $M^n$ and by $\widetilde{\nabla}$ that
of $\M^{n+k}$. With the appropriate identifications, the formulas of
Gauss and Weingarten state respectively
\begin{eqnarray}
\widetilde{\nabla}_X Y &=& \nabla_X Y + \alpha(X,Y), \label{3.1} \\
\widetilde{\nabla}_X \xi &=& -S_{\xi}X + \nabla^{\perp}_X\xi,
\label{3.2}
\end{eqnarray}
where $X$ and $Y$ are vector fields tangent to $M^n$ and $\xi$ is a
normal vector field along $M^n$. The symmetric (1,2)-tensor field
$\alpha$, taking values in the normal bundle, is called the
\emph{second fundamental form}, the symmetric (1,1)-tensor field
$S_{\xi}$ on $M^n$ is the \emph{shape operator associated to $\xi$}
and $\nabla^{\perp}$ is a connection in the normal bundle. From
these formulas the equations of Gauss and Codazzi can be deduced:
\begin{eqnarray}
\mathrm{tan}(\widetilde{R}(X,Y)Z) &=& R(X,Y)Z+S_{\alpha(X,Z)}Y-S_{\alpha(Y,Z)}X, \label{3.3} \\
\mathrm{tan}(\widetilde{R}(X,Y)\xi) &=&
(\nabla_YS)_{\xi}X-(\nabla_XS)_{\xi}Y, \label{3.4}
\end{eqnarray}
for $p\in M^n$ and $X,Y,Z\in T_pM^n$, $\xi\in T_p^{\perp}M^n$. Here
$R$ is the Riemann-Christoffel curvature tensor of $M^n$,
$\widetilde{R}$ that of $\M^{n+k}$, ``tan'' denotes the projection
on the tangent space to $M^n$ and
$(\nabla_XS)_{\xi}Y=\nabla_X(S_{\xi}Y)-S_{\xi}(\nabla_XY)-S_{\nabla_X^{\perp}\xi}Y$.\medskip

Now let $F:M^2\rightarrow\M^3(\kappa,\tau)$ be an isometric
immersion of an oriented surface in a BCV space, with unit normal
$\xi$ and associated shape operator $S$. We denote by $\theta$ the
angle between $e_3$ and $\xi$ and by $T$ the projection of $e_3$ on
the tangent plane to $M^2$, i.e. the vector field $T$ on $M^2$ such
that $F_{\ast}T+\cos\theta\,\xi=e_3$. If we work locally, we may
assume $\theta\in [0,\frac{\pi}{2}]$. The equations of Gauss
(\ref{3.3}) and Codazzi (\ref{3.4}) give respectively
\begin{multline}\label{Gaussvgl}
R(X,Y)Z = (\kappa-3\tau^2)(\langle Y,Z\rangle X-\langle X,Z\rangle
          Y)
          -(\kappa-4\tau^2)(\langle Y,T \rangle\langle Z,T \rangle X - \langle X,T \rangle\langle Z,T \rangle Y\\
          +\langle X,T \rangle\langle Y,Z \rangle T - \langle Y,T \rangle\langle X,Z \rangle T)
          +\langle SY,Z \rangle SX - \langle SX,Z \rangle SY
\end{multline}
and
\begin{equation}
\nabla_XSY-\nabla_YSX-S[X,Y]=(\kappa-4\tau^2)\cos\theta(\langle
Y,T\rangle X-\langle X,T\rangle Y)\label{Codazzivgl}
\end{equation}
for $p\in M^2$ and $X,Y,Z\in T_pM^2$. From (\ref{Gaussvgl}) it
follows moreover that the Gaussian curvature of $M^2$ is given by
\begin{equation}
K=\det S+\tau^2+(\kappa-4\tau^2)\cos^2\theta.\label{Gaussvgl2}
\end{equation}
Finally, we remark that the following structure equations hold for
$p\in M^2$ and $X\in T_pM^2$:
\begin{eqnarray}
\nabla_XT &=& \cos\theta(SX-\tau JX),\label{structuurvgl1}\\
X[\cos\theta] &=& -\langle SX-\tau JX,T\rangle,\label{structuurvgl2}
\end{eqnarray}
where $J$ denotes the rotation over $\frac{\pi}{2}$ in $T_pM^2$.
These equations can be verified straightforwardly by comparing the
tangential and normal components of both sides of the equality
$\nt_X(T+\cos\theta\,\xi)= \tau(X\times(T+\cos\theta\,\xi))$.\\

The following theorem is proven in \cite{7}:
\begin{theorem}\label{stelling Daniel}\cite{7}
Let $M^2$ be a simply connected, oriented Riemannian surface with
metric $\langle\cdot,\cdot\rangle$, Levi Civita connection $\nabla$
and curvature tensor $R$. Let $J$ denote the rotation over
$\frac{\pi}{2}$ in $TM^2$ and $S$ a field of symmetric operators on
$TM^2$. Finally, let $T$ be a vector field on $M^2$ and let
$\cos\theta$ be a differentiable function, satisfying $\langle
T,T\rangle + \cos^2\theta =1$. Then there exists an isometric
immersion $F$ of $M^2$ in $\M^3(\kappa,\tau)$ with unit normal
$\xi$, such that $S$ is the shape operator and
$e_3=F_{\ast}T+\cos\theta\,\xi$ if and only if the equations
(\ref{Gaussvgl}), (\ref{Codazzivgl}), (\ref{structuurvgl1}) and
(\ref{structuurvgl2}) are satisfied. In this case the immersion is
moreover unique up to a global isometry of $\M^3(\kappa,\tau)$,
preserving both the orientations of the base space $\M^2(\kappa)$
and the fibres of $\pi$.
\end{theorem}
\medskip

\section{Parallel, semi-parallel and higher order parallel hypersurfaces}

Let $F:M^n\rightarrow\widetilde{M}^{n+1}$ be an isometric immersion
of Riemannian manifolds and $p\in M^n$. If $\alpha$ is the second
fundamental form and $\xi$ is a unit normal vector field on the
hypersurface, we define the scalar valued second fundamental form
$h$ to be the (0,2)-tensor field satifying $\alpha(X,Y)=h(X,Y)\,\xi$
for all $p\in M^n$ and $X,Y\in T_pM^n$. The covariant derivative of
$h$ is defined by
$$(\nabla h)(X,Y,Z) = X[h(Y,Z)] - h(\nabla_XY,Z) - h(Y,\nabla_XZ),$$
for all $X,Y,Z\in T_pM^n$ with $\nabla$ the Levi Civita connection
of $M^n$. If $R$ is the curvature tensor of $M^n$, we also define
$$(R\cdot h)(X,Y,Z_1,Z_2)=-h(R(X,Y)Z_1,Z_2)-h(Z_1,R(X,Y)Z_2),$$
for all $X,Y,Z_1,Z_2\in T_pM^n$. If $\nabla h=0$, we say that $M^n$
has parallel second fundamental form or, for short, that it is a
\emph{parallel} hypersurface. If $R\cdot h=0$, we say that $M^n$ is
a \emph{semi-parallel} hypersurface.\\

For any integer $k\geq2$, we define recursively
\begin{multline*}(\nabla^kh)(X_1,\ldots,X_k,Y,Z) = X_1[(\nabla^{k-1}h)(X_2,\ldots,X_k,Y,Z)]\\
-(\nabla^{k-1}h)(\nabla_{X_1}X_2,\ldots,X_k,Y,Z) - \ldots
-(\nabla^{k-1}h)(X_2,\ldots,X_k,Y,\nabla_{X_1}Z)\end{multline*} for
$X_1,\ldots, X_k,Y,Z\in T_pM^n$. We call a hypersurface satisfying
$\nabla^kh=0$ a \emph{$k$-parallel} hypersurface or a \emph{higher
order parallel} hypersurface. With slight modifications, all these
notions can also be defined for
submanifolds with arbitrary codimension.\\

The classification of parallel hypersurfaces in real space forms is
proven in \cite{14}, whereas for the classification of $k$-parallel
hypersurfaces in real space forms we refer to \cite{8}, \cite{9} and
\cite{10}:
\begin{theorem} \cite{14} A parallel hypersurface in a simply connected, complete real
space form of constant sectional curvature $c$ is one of the
following. In $\E^{n+1}$: an open part of a product immersion
$\E^k\times\S^{n-k}$, $k\in\{0,\ldots, n\}$. In $\S^{n+1}(c)$: an
open part of a product immersion $\S^k\times\S^{n-k}$,
$k\in\{0,\ldots, n\}$. In $\H^{n+1}(c)$: an open part of a product
immersion $\H^k\times\S^{n-k}$, $k\in\{0,\ldots, n\}$ or of a
horosphere.
\end{theorem}
\begin{theorem} \cite{8}, \cite{9}, \cite{10}
A $k$-parallel hypersurface in a simply connected, complete real
space form of constant sectional curvature $c$ is one of the
following. In $\E^{n+1}$: an open part of a parallel hypersurface or
of a cylinder on a plane curve, whose curvature is a polynomial
function of degree at most $k-1$ of the arc length. In
$\S^{n+1}(c)$: an open part of a parallel hypersurface or, for
$n=2$, of the inverse image under the Hopf-fibration
$\S^3(c)\rightarrow\S^2(4c)$ of a spherical curve in $\S^2(4c)$
whose geodesic curvature is a polynomial of degree at most $k-1$ of
the arc length. In $\H^{n+1}(c)$: an open part of a parallel
hypersurface.
\end{theorem}
\medskip

In \cite{3} the following classification for parallel surfaces in
BCV spaces is proven:
\begin{theorem} \label{theo4} \cite{3}
A parallel surface in $\M^3(\kappa,\tau)$, with $\kappa\neq
4\tau^2$, is an open part of a Hopf cylinder over a Riemannian
circle in $\M^2(\kappa)$ or of a totally geodesic leaf of the Hopf
fibration, the latter case only occuring for $\tau=0$.
\end{theorem}
The technique used in the proof of this theorem is based on the fact
that for parallel surfaces the left-hand side of Codazzi's equation
(\ref{Codazzivgl}) is zero. For $k$-parallel
surfaces another approach is needed.\\

We refer to \cite{11} for a proof of the following lemma:
\begin{lemma}\label{lemma1} \cite{11}
A $k$-parallel surface immersed in a three-dimensional Riemannian
manifold is semi-parallel, or equivalently, it is flat or totally
umbilical.
\end{lemma}
This means that in our search for $k$-parallel surfaces in BCV
spaces, we can focus on totally umbilical surfaces (meaning that at
every point the shape operator is a scalar multiple of the identity)
and flat surfaces (meaning that the Gaussian curvature at every
point is zero). In the next section we will give a complete
classification of totally umbilical surfaces in BCV spaces and in
the last section we will classify all flat, $k$-parallel surfaces in
BCV spaces.

\section{Totally umbilical surfaces}

In \cite{17}, it was proven that there are no totally umbilical
surfaces in the Heisenberg group $\Nil$. The following lemma
generalizes this result.
\begin{lemma}\label{lemma2}
Let $M^2\rightarrow\M^3(\kappa,\tau)$ be a totally umbilical surface
with shape operator $S=\lambda\,\mathrm{id}$. Then $\tau=0$ and the
following equations hold:
\begin{equation}\label{5.1}
T[\lambda]=-\kappa\cos\theta\sin^2\theta, \quad (JT)[\lambda]=0,
\quad T[\theta]=\lambda\sin\theta, \quad (JT)[\theta]=0,
\end{equation}
\begin{equation}\label{5.2}
\nabla_T T=\lambda\cos\theta\,T, \quad
\nabla_{JT}T=\lambda\cos\theta\,JT, \quad \nabla_T
JT=\lambda\cos\theta\,JT, \quad \nabla_{JT}JT=-\lambda\cos\theta\,T.
\end{equation}
\end{lemma}

\noindent\emph{Proof.} First assume that $\theta$ is identically
zero. Then with the notations of section 2 we have
$TM^2=\mathrm{span}\{e_1,e_2\}$. But according to Frobenius' theorem
and (\ref{Lie-haken}), this distribution is only integrable if
$\tau=0$. Now $T=JT=0$ and, since $Se_1=-\nt_{e_1}e_3=0$ and
$Se_2=-\nt_{e_2}e_3=0$, also $\lambda=0$. All equations stated in
the lemma are satisfied.

We now work on an open subset of $M^2$ where $\theta$ is nowhere
zero. From Codazzi's equation (\ref{Codazzivgl}) for $X=T$ and
$Y=JT$, we get
\begin{equation}\label{5.3}
T[\lambda]=-(\kappa-4\tau^2)\cos\theta\sin^2\theta, \quad
JT[\lambda]=0.
\end{equation}
The structure equations (\ref{structuurvgl1}) and
(\ref{structuurvgl2}) yield
\begin{equation}\label{5.4}
\nabla_T T=\cos\theta(\lambda T - \tau JT), \quad
\nabla_{JT}T=\cos\theta(\tau T + \lambda JT), \quad
T[\theta]=\lambda\sin\theta, \quad (JT)[\theta]=\tau\sin\theta.
\end{equation}
Using orthonormal expansion and $\langle T,JT\rangle=0$, $\langle
T,T\rangle = \langle JT,JT\rangle = \sin^2\theta$, we get
\begin{equation}\label{5.5}
\nabla_T JT = \cos\theta(\tau T + \lambda JT), \quad
\nabla_{JT}JT=\cos\theta(-\lambda T + \tau JT).
\end{equation}

Remark that $[T,JT]=\nabla_T JT - \nabla_{JT} T = 0$ and hence
$$0=[T,JT][\lambda]=T[(JT)[\lambda]]-(JT)[T[\lambda]]=(\kappa-4\tau^2)\tau\sin^2\theta(2\cos^2\theta-\sin^2\theta).$$
Since we assume $\kappa-4\tau^2\neq 0$ and $\sin\theta\neq 0$,
either $\tau=0$ or $2\cos^2\theta-\sin^2\theta=0$. But the latter
implies that $\theta$ is a constant and then from the last equation
of (\ref{5.4}) we also get $\tau=0$. The equations stated in the
lemma follow easily from (\ref{5.3}), (\ref{5.4}) and (\ref{5.5}).
By a continuity argument, these will hold on the whole of $M^2$.
\hfill$\square$\bigskip

The following is an immediate corollary of Lemma \ref{lemma2}.

\begin{theorem} \label{theo} The only BCV spaces admitting totally umbilical
surfaces are the Riemannian products
$(\S^2(\kappa)\setminus\{\infty\})\times\R$ and
$\H^2(\kappa)\times\R$.
\end{theorem}

It is now sufficient to study totally umbilical surfaces in
$\S^2(\kappa)\times\R$ and $\H^2(\kappa)\times\R$. To do this, we
consider these spaces as hypersurfaces of the four-dimensional
Euclidean space $\E^4$ and the four-dimensional Lorentzian space
$\L^4=(\R^4,-dx_1^2+dx_2^2+dx_3^2+dx_4^2)$ respectively:
$$\S^2(\kappa)\times\R=\left\{(x_1,x_2,x_3,x_4)\in\E^4\ \left|\ x_1^2+x_2^2+x_3^2=\frac{1}{\kappa}\right.\right\}$$
and
$$\H^2(\kappa)\times\R=\left\{(x_1,x_2,x_3,x_4)\in\L^4\ \left|\ -x_1^2+x_2^2+x_3^2=\frac{1}{\kappa},\ x_1>0\right.\right\}.$$
Remark that in both cases the vector field $\widetilde{\xi}$,
defined by
$\widetilde{\xi}(x_1,x_2,x_3,x_4)=\sqrt{|\kappa|}(x_1,x_2,x_3,0)$,
is orthogonal to the hypersurface and that
$\langle\widetilde{\xi},\widetilde{\xi}\rangle=1$ in the first case
and $\langle\widetilde{\xi},\widetilde{\xi}\rangle=-1$ in the second
case.\\

First, we remark that the only totally umbilical surfaces that are
also higher order parallel are trivial:

\begin{proposition} \label{prop1}
A $k$-parallel, totally umbilical surface in $\S^2(\kappa)\times\R$,
respectively $\H^2(\kappa)\times\R$, is totally geodesic and an open
part of $\S^2(\kappa)\times\{t_0\}$ or $\S^1(\kappa)\times\R$,
respectively $\H^2(\kappa)\times\{t_0\}$ or $\H^1(\kappa)\times\R$.
Moreover, these surfaces are the only totally geodesic ones.
\end{proposition}

\noindent\emph{Proof.} If $\theta$ is identically zero, the surface
is an open part of $\S^2(\kappa)\times\{t_0\}$ or
$\H^2(\kappa)\times\{t_0\}$. Hence we may assume that $\theta\neq
0$. Putting $U=\frac{T}{\|T\|}=\frac{T}{\sin\theta}$ and $V=JT$, we
have $[U,V]=0$, so we can take coordinates $(u,v)$ with
$U=\frac{\partial}{\partial u}$ and $V=\frac{\partial}{\partial v}$.
Remark that $\lambda$ and $\theta$ only depend on $u$ and
\begin{equation} \label{5.6}
\lambda'=-\kappa\cos\theta\sin\theta=-\frac{\kappa}{2}\sin(2\theta),
\quad \theta'=\lambda.
\end{equation}
Since
$$\nabla_U U=\frac{1}{\sin\theta}\left(T\left[\frac{1}{\sin\theta}\right]T+\frac{1}{\sin\theta}\nabla_T T\right)=0,$$
we have
$$0=(\nabla^kh)(U,U,\ldots,U,U)=U[U[\ldots U[h(U,U)]\ldots]]=\lambda^{(k)}(u),$$
which implies that $\lambda$ is a polynomial of degree at most $k-1$
in $u$. Now from (\ref{5.6}), we see that both $\sin(2\theta)$ and
$\theta$ are polynomials in $u$. The only possibility is that
$\theta$ is a constant and thus, again from (\ref{5.6}), $\lambda=0$
and $\cos\theta=0$. So $\theta=\frac{\pi}{2}$ and the surface is an
open part of $\gamma\times\R$, with $\gamma$ a curve in
$\S^2(\kappa)$ or $\H^2(\kappa)$.

It remains to prove that $\gamma$ is a geodesic. We continue the
proof for $\kappa>0$, but the other case is completely similar.
Assume that $\gamma$ is parametrized by arc length and denote the
immersion by
$$F:M^2\rightarrow\S^2(\kappa)\times\R\subset\E^4:(s,t)\mapsto(\gamma(s),t).$$
Denoting by `` $\cdot$ '' the inner product on $\E^3$ and by ``
$\times$ '' the cross product, we have that $F_s=(\gamma',0)$ and
$F_t=(0,1)$ span the tangent space, that
$\widetilde{\xi}=\sqrt{\kappa}(\gamma,0)$ is a unit vector
orthogonal to the surface and orthogonal to $\S^2(\kappa)\times\R$
and $\xi=\sqrt{\kappa}(\gamma\times\gamma',0)$ is a unit vector
orthogonal to $M^2$, tangent to $\S^2(\kappa)\times\R$. Moreover
\begin{eqnarray*}
\left\langle S\frac{\partial}{\partial s},\frac{\partial}{\partial s}\right\rangle &=& \langle F_{ss},\xi\rangle = \kappa((\gamma\times\gamma')\cdot\gamma'',0),\\
\left\langle S\frac{\partial}{\partial s},\frac{\partial}{\partial t}\right\rangle &=& \langle F_{st},\xi\rangle = (0,0),\\
\left\langle S\frac{\partial}{\partial t},\frac{\partial}{\partial
t}\right\rangle &=& \langle F_{tt},\xi\rangle = (0,0),
\end{eqnarray*}
and thus the surface is totally umbilical (and automatically totally
geodesic) if and only if $(\gamma\times\gamma')\cdot\gamma''=0$, or
equivalently, if and only if $\gamma''$ is proportional to $\gamma$.
This means that $\gamma''$ has no component tangent to
$\S^2(\kappa)$ and hence has to be a geodesic, i.e. a great circle.

The fact that these surfaces are the only totally geodesic ones
follows immediately from the first equation of (\ref{5.1}). \hfill
$\square$ \medskip

Before proceeding with the full classification, we develop some
machinery to study surfaces in $\S^2(\kappa)\times\R$ and
$\H^2(\kappa)\times\R$.\medskip

Consider an isometric immersion
$F:M^2\rightarrow\S^2(\kappa)\times\R$. Denoting by $\xi$ a unit
vector tangent to $\S^2(\kappa)\times\R$ and normal to $M^2$, one
easily sees that the fourth components of $F_{\ast}T$, $F_{\ast}JT$
and $\xi$ in $\E^4$ satisfy
\begin{equation}\label{5.6a}
(F_{\ast}T)_4=\sin^2\theta, \quad (F_{\ast}JT)_4=0, \quad
\xi_4=\cos\theta.
\end{equation}
Take $\widetilde{\xi}$ as above and let $X$ be a tangent vector to
$M^2$. Then $\langle
\nabla^{\perp}_X\widetilde{\xi},\xi\rangle=\langle
D_X\widetilde{\xi},\xi\rangle=X_1\xi_1+X_2\xi_2+X_3\xi_3=-X_4\xi_4=-\langle
X,T\rangle\cos\theta,$ where $D$ denotes the Euclidean connection.
Thus, the normal connection of $M^2$ as a submanifold of $\E^4$ is
given by
\begin{equation*}
\nabla^{\perp}_X\widetilde{\xi} = -\langle
X,T\rangle\cos\theta\,\xi, \quad \nabla^{\perp}_X\xi = \langle
X,T\rangle\cos\theta\,\widetilde{\xi}.
\end{equation*}
Using Weingarten's formula, we see that the shape operator
associated to $\widetilde{\xi}$, which we denote by $\widetilde{S}$,
must satisfy
\begin{eqnarray*}
F_{\ast}(\widetilde{S}T) &=& (-(F_{\ast}T)_1,-(F_{\ast}T)_2,-(F_{\ast}T)_3,0)-\cos\theta\sin^2\theta\,(\xi_1,\xi_2,\xi_3,\cos\theta),\\
F_{\ast}(\widetilde{S}(JT)) &=&
(-(F_{\ast}JT)_1,-(F_{\ast}JT)_2,-(F_{\ast}JT)_3,0)=-JT.
\end{eqnarray*}
The second equation implies that the matrix of $\widetilde{S}$ with
respect to the basis $\{T,JT\}$ takes the form
$$\widetilde{S}=\left(\begin{array}{cc} a&0 \\ 0&-1 \end{array}\right)$$
and looking at the fourth component of the first equation we get
$a=-\cos^2\theta$. Hence
\begin{equation} \label{5.7}
\widetilde{S}=\left(\begin{array}{cc} -\cos^2\theta&0 \\ 0&-1
\end{array}\right).
\end{equation}
Remark that from the other components of the first equation
\begin{equation} \label{5.8}
(F_{\ast}T)_j=-\cos\theta\,\xi_j,\qquad j=1,2,3.
\end{equation}
\medskip

We can do exactly the same for $\H^2(\kappa)\times\R$. The equations
(\ref{5.6a}) remain the same. The normal connection changes to
\begin{equation*}
\nabla^{\perp}_X\widetilde{\xi}=-\langle X,T\rangle\cos\theta\,\xi,
\quad \nabla^{\perp}_X\xi = -\langle
X,T\rangle\cos\theta\,\widetilde{\xi},
\end{equation*}
but the shape operator associated to $\widetilde{\xi}$, (\ref{5.7}),
and formula (\ref{5.8}) remain the same.\\

We will now classify totally umbilical surfaces in $\S^2(1)\times\R$
and $\H^2(-1)\times\R$ and for arbitrary $\kappa$ the totally
umbilical surfaces will then be homothetic to these.

\begin{theorem} \label{T:3.3.2}
Let $F:M^2\rightarrow\S^2(1)\times\R\subset\E^4$ be a totally
umbilical surface with shape operator $S=\lambda\,\mathrm{id}$ and
angle function $\theta$, which is not totally geodesic. Then one can
choose local coordinates $(u,v)$ on $M^2$ such that $\lambda$ and
$\theta$ only depend on $u$ and
\begin{equation}\label{3.3.13}\theta(u)=\arctan\left(\frac{2ce^{\pm cu}}{1-c^2+e^{\pm 2cu}}\right), \qquad \lambda(u)=\frac{\theta'(u)}{\sin\theta(u)},\end{equation}
for some real constant $c>0$. Moreover, the immersion is, up to an
isometry, locally given by
\begin{equation}\label{3.3.14}
F(u,v)=\frac{1}{c}\left(\lambda,\,\sin\theta\,\cos
v,\,\sin\theta\,\sin v,\,c\int\sin^2\theta\,du\right).
\end{equation}
\end{theorem}

\noindent\emph{Proof.} It follows from \e{5.2} that $[T,JT]=0$.
Hence, we can take local coordinates $(u,v)$ on $M^2$, such that
$T=\frac{\partial}{\partial u}$, $JT=\frac{\partial}{\partial v}$.
From \e{5.1} we see that $\lambda$ and $\theta$ only depend on $u$
and that they satisfy
\begin{equation}\label{3.3.15}
\lambda^2+\sin^2\theta=c^2, \qquad \theta'=\lambda\sin\theta,
\end{equation}
for some strictly positive real constant $c$.

From the formula of Gauss, \e{5.1}, \e{5.7} and \e{5.8}, we obtain
for $j=1,2,3$
\begin{eqnarray}
(F_j)_{uu} &=& \lambda\cos\theta(F_j)_u-\lambda\frac{\sin^2\theta}{\cos\theta}(F_j)_u-\cos^2\theta\sin^2\theta\,F_j, \label{3.3.16}\\
(F_j)_{uv} &=& \lambda\cos\theta(F_j)_v, \label{3.3.17}\\
(F_j)_{vv} &=& -\frac{\lambda}{\cos\theta}(F_j)_u-\sin^2\theta\,F_j.
\label{3.3.18}
\end{eqnarray}
The equations for the fourth component are trivially satisfied. The
solution of (\ref{3.3.17}) is
\begin{equation}\label{3.3.19}
F_j=(A_j(u)+B_j(v))\exp\left(\int\lambda\cos\theta\,du\right),
\end{equation}
where $A_j$ and $B_j$ are real-valued functions in one variable.
Substituting this in (\ref{3.3.16}) yields
\begin{equation}\label{3.3.20}
A_j=a_j\int\exp\left(-\int\frac{\lambda}{\cos\theta}\,du\right)\,du
+ \alpha_j
\end{equation}
with $a_j$, $\alpha_j\in\R$, and substituting it in (\ref{3.3.18})
gives $B_j''+(\lambda^2+\sin^2\theta)B_j=
A_j''-(\lambda^2+\sin^2\theta)A_j$, or equivalently $B_j''+c^2B_j=
A_j''-c^2A_j$. It is easy to check that the right hand side of this
equation is constant and thus the solution for $B_j$ is
\begin{equation}\label{3.3.21}
B_j=b_j\cos(cv)+\beta_j\sin(cv)+\frac{A_j''}{c^2}-A_j,
\end{equation}
with $b_j$, $\beta_j\in\R$. By substituting \e{3.3.20} and
\e{3.3.21} in \e{3.3.19}, we conclude that the functions $F_j$ take
the form
\begin{multline} \label{3.3.22}
F_j=\left(-a_j\frac{\lambda}{c^2\cos\theta}\exp\left(-\int\frac{\lambda}{\cos\theta}\,du
\right)+b_j\cos(cv)\right.\\
+\beta_j\sin(cv)\left)\,
\exp\left(\int\lambda\cos\theta\,du\right)\right., \quad j=1,2,3
\end{multline}
and from (\ref{5.6a}):
\begin{equation} \label{3.3.23}
F_4=\int\sin^2\theta\,du.
\end{equation}

There are some conditions on $F$ which we have neglected so far,
namely $F\in\S^2(1)\times\R$,  $\langle \xi,F_u\rangle=\langle
\xi,F_v \rangle=0$,  $\langle \widetilde{\xi},F_u \rangle=\langle
\widetilde{\xi},F_v \rangle=0$, $\langle F_u,F_u \rangle = \langle
F_v,F_v \rangle = \sin^2\theta$, $\langle \xi,\xi \rangle = \langle
\widetilde{\xi},\widetilde{\xi} \rangle =1$ and $\langle
\xi,\widetilde{\xi} \rangle = \langle F_u,F_v \rangle = 0$. These
are equivalent to
\begin{equation}  \label{3.3.24}
\sum_{j=1}^3 F_j^2 = 1, \quad \sum_{j=1}^3 (F_j)_u^2 =
\cos^2\theta\sin^2\theta, \quad \sum_{j=1}^3 (F_j)_v^2 =
\sin^2\theta, \quad \sum_{j=1}^3 (F_j)_u(F_j)_v = 0.
 \end{equation} Now looking at $a=(a_1,a_2,a_3)$,
$b=(b_1,b_2,b_3)$ and $\beta=(\beta_1,\beta_2,\beta_3)$ as vectors
in $\R^3$ with the Euclidean inner product `` $\cdot$ '', the
conditions (\ref{3.3.24}) are equivalent to
\begin{equation*}
a\cdot b = a\cdot\beta = b\cdot\beta=0,
\end{equation*}
\begin{equation*}
\|a\|^2=a\cdot a=
c^2\cos^2\theta\,\exp\left(2\int\frac{\lambda\sin^2\theta}{\cos\theta}\,du\right),
\end{equation*}
\begin{equation*}
\|b\|^2=b\cdot b=\beta\cdot\beta
=\frac{\sin^2\theta}{c^2}\,\exp\left(-2\int\lambda\cos\theta\,du\right).
\end{equation*}
Remark that the right hand sides of these equations are constant.
They imply that after a suitable isometry of $\S^2(1)\times\R$ we
may assume that
\begin{equation*}
a =
\left(-c\cos\theta\exp\left(\int\frac{\lambda\sin^2\theta}{\cos\theta}\,du\right),\,0,\,0\right),
\end{equation*}
\begin{equation*}
b =
\left(0,\,\frac{\sin\theta}{c}\exp\left(-\int\lambda\cos\theta\,du\right),\,0\right),
\end{equation*}
\begin{equation*}
\beta
=\left(0,\,0,\,\frac{\sin\theta}{c}\exp\left(-\int\lambda\cos\theta\,du\right)\right).
\end{equation*}
Now the reparametrization $cv\mapsto v$ gives the result \e{3.3.14}.

To conclude, we solve the equations \e{3.3.15} explicitly. Putting
$\theta=\arctan(f)$, we obtain
$$\left(\frac{\theta'}{\sin\theta}\right)^2+\sin^2\theta=c^2
\Leftrightarrow
\left(\frac{f'}{f\sqrt{1+f^2}}\right)^2+\frac{f^2}{1+f^2}=c^2
\Leftrightarrow \frac{(f')^2}{f^2(c^2+(c^2-1)f^2)}=1.$$ From the
last equation we see that $c^2+(c^2-1)f^2$ has to be positive and
hence we can proceed by integration:
\begin{eqnarray*}
\frac{f'}{f\sqrt{c^2+(c^2-1)f^2}}=\pm 1 &\Leftrightarrow& \ln\left(\frac{c+\sqrt{c^2+(c^2-1)f^2}}{f}\right)=\pm cu+d\\
&\Leftrightarrow& f=\frac{2c\,e^{\pm cu+d}}{1-c^2+e^{2(\pm cu+d)}},
\end{eqnarray*}
for some $d\in\R$. After a change of the $u$-coordinate, which does
not change $\frac{\partial}{\partial u}$, we obtain the result
\e{3.3.13}. \hfill $\square$\medskip

\begin{remark}\label{R:3.3.1} We can write \e{3.3.14} in a more explicit form. After the reparametrization $e^{\pm cu}\mapsto u$ and,
if necessary, an isometry switching the sign of some of the
components, (\ref{3.3.14}) is given by
\begin{equation*}
F(u,v)=\left(\frac{2u\cos v}{p(u)q(u)},\ \frac{2u\sin v}{p(u)q(u)},\
\frac{1-c^2-u^2}{p(u)q(u)},\
\ln\left(\frac{p(u)}{q(u)}\right)\right),
\end{equation*}
where $p(u)=\sqrt{u^2+(c-1)^2}$ and $q(u)=\sqrt{u^2+(c+1)^2}$.
\end{remark}\medskip

\begin{theorem}\label{T:3.3.3}
Let $F:M^2\rightarrow\H^2(-1)\times\R\subset\E^4_1$ be a totally
umbilical surface with shape operator $S=\lambda\,\mathrm{id}$ and
angle function $\theta$, which is not totally geodesic. Then one can
choose local coordinates $(u,v)$ on $M^2$ such that $\lambda$ and
$\theta$ only depend on $u$ and we are in one of the following three
cases:
\begin{itemize}
\item[$\mathrm{(i)}$] $\theta(u)$ and $\lambda(u)$ are given by
\begin{equation}\label{3.3.26}\theta(u)=\arctan\left(\frac{2ce^{\pm cu}}{1+c^2-e^{\pm 2cu}}\right),
\qquad \lambda(u)=\frac{\theta'(u)}{\sin\theta(u)},\end{equation}
for some real constant $c>0$, and the immersion is, up to an
isometry, locally given by
\begin{equation} \label{3.3.27}
F(u,v)=\frac{1}{c}\left(\lambda,\,\sin\theta\,\cos
v,\,\sin\theta\,\sin v,\,c\int\sin^2\theta\,du\right),
\end{equation}
\item[$\mathrm{(ii)}$] $\theta(u)$ and $\lambda(u)$ are given by
\begin{equation}\label{3.3.28}\theta(u)=\mathrm{arccot}\,(\pm u),
\qquad \lambda(u)=\frac{\mp 1}{\sqrt{1+u^2}},\end{equation} and the
immersion is, up to an isometry, locally given by
\begin{equation} \label{3.3.29}
F(u,v)=\frac{1}{\sqrt{1+u^2}}\left(\frac{u^2+v^2}{2}+1,\,v,\,\frac{u^2+v^2}{2},\,\sqrt{1+u^2}\arctan
u\right),
\end{equation}
\item[$\mathrm{(iii)}$] $\theta(u)$ and $\lambda(u)$ are given by
\begin{equation}\label{3.3.30}\theta(u)=\arctan\left(\frac{\tan c}{\sin(\pm u\sin c)}\right),
\qquad \lambda(u)=\frac{\theta'(u)}{\sin\theta(u)},\end{equation}
for some real constant $c\neq 0$, and the immersion is, up to an
isometry, locally given by
\begin{equation}\label{3.3.31}
F(u,v)=\frac{1}{\sin c}\left(\sin\theta\,\cosh v,\,\sin\theta\,\sinh
v,\,\lambda,\,\sin c\int\sin^2\theta\,du\right).
\end{equation}
\end{itemize}
\end{theorem}

\noindent\emph{Proof.} We use again the coordinates $(u,v)$ such
that $T=\frac{\partial}{\partial u}u$ and
$JT=\frac{\partial}{\partial v}$. From \e{5.1}, we obtain that
$\lambda$ and $\theta$ only depend on $u$ and that they satisfy
$\lambda^2-\sin^2\theta=C$, $\theta'=\lambda\sin\theta$, for some
real constant $C>-1$. The formula of Gauss yields for $j=1,2,3$:
\begin{eqnarray*}
(F_j)_{uu} &=& \lambda\cos\theta(F_j)_u-\lambda\frac{\sin^2\theta}{\cos\theta}(F_j)_u+\cos^2\theta\sin^2\theta F_j,\\
(F_j)_{uv} &=& \lambda\cos\theta(F_j)_v,\\
(F_j)_{vv} &=& -\frac{\lambda}{\cos\theta}(F_j)_u+\sin^2\theta F_j,
\end{eqnarray*}
such that $F_j$ again takes the form (\ref{3.3.19}), with $A_j$
again equal to (\ref{3.3.20}). The differential equation for $B_j$
becomes $B_j''+(\lambda^2-\sin^2\theta)B_j=
A_j''-(\lambda^2-\sin^2\theta)A_j$, or equivalently
$B_j''+CB_j=A_j''-CA_j$. The right hand side of this equation is
again constant. We now consider three cases.

\vskip.05in\emph{Case} (A): $C>0$. This case corresponds to the
first case of the theorem. We can put $C=c^2$ for some strictly
positive real constant $c$. The rest of the proof is similar to the
one above and we will therefore omit it.

\vskip.05in\emph{Case} (B): $C=0$. The solution of the equations
$\lambda^2=\sin^2\theta$ and $\theta'=\lambda\sin\theta$ is given by
\e{3.3.28}. Substituting this in \e{3.3.20} yields that $A_j$ takes
the form $A_j(u)=p_ju^2+q_j$ for some $p_j,q_j\in\R$. The equation
for $B_j$ becomes $B_j''=A_j''$. From this equation and \e{3.3.19},
we obtain
$$F_j=\frac{1}{\sqrt{1+u^2}}(a_j(u^2+v^2)+b_jv+c_j), \qquad
j=1,2,3,$$ where $a_j,b_j,c_j\in\R$. Moreover, from \e{5.6a} and
\e{3.3.28}, we have
$$F_4=\int\sin^2\theta\,du=\arctan u.$$

The conditions analogous to (\ref{3.3.24}) now read
\begin{equation}\label{3.3.32}
\begin{aligned}
& -F_1^2+F_2^2+F_3^2=-1, \\
& -(F_1)_u^2+(F_2)_u^2+(F_3)_u^2=\cos^2\theta\sin^2\theta,\\
& -(F_1)_v^2+(F_2)_v^2+(F_3)_v^2=\sin^2\theta, \\
& -(F_1)_u(F_1)_v+(F_2)_u(F_2)_v+(F_3)_u(F_3)_v=0,
\end{aligned}
\end{equation}
and looking at $a$, $b$ and $c$ as vectors in $\R^3$, but now
equipped with the standard Lorentzian inner product `` $\cdot$ '',
these are equivalent to $a\cdot a=a\cdot b=b\cdot c=0$, $a\cdot
c=-\frac{1}{2}$, $b\cdot b=1$, $c\cdot c=-1$. After a suitable
isometry of $\H^2(-1)\times\R$, we may assume that
$a=(\frac{1}{2},0,\frac{1}{2})$, $b=(0,1,0)$ and $c=(1,0,0)$. This
gives the result \e{3.3.29}.

\vskip.05in\emph{Case} (C): $C<0$. Clearly, we have $C>-1$ and hence
we may put $C=-\sin^2c$, for some real number $c$. The equation for
$B_j$ becomes $B_j''-\sin^2c\,B_j=A_j''+\sin^2c\,A_j$, with solution
\begin{equation}\label{3.3.33}
B_j=b_j\cosh(v\sin c)+\beta_j\sinh(v\sin
c)-\frac{A_j''}{\sin^2c}-A_j.
\end{equation}
Hence $F$ is given by
\begin{multline} \label{3.3.34}
F_j=\left(-a_j\frac{\lambda}{(\sin^2c)\cos\theta}\exp\left(-\int\frac{\lambda}{\cos\theta}\,du
\right)+b_j\cosh(v\sin c)\right.\\
+\beta_j\sinh(v\sin c)\left)\,
\exp\left(\int\lambda\cos\theta\,du\right)\right., \quad j=1,2,3
\end{multline}
and $F_4$ takes the form \e{3.3.23}.

Looking at $a$, $b$ and $\beta$ as vectors in $\R^3$ with the
standard Lorentzian inner product, the conditions \e{3.3.32} yield
\begin{equation*}
a\cdot b = a\cdot\beta = b\cdot\beta=0,
\end{equation*}
\begin{equation*}
a\cdot a
=\sin^2c\,\cos^2\theta\,\exp\left(2\int\frac{\lambda\sin^2\theta}{\cos\theta}\,du\right),
\end{equation*}
\begin{equation*}
b\cdot b=-\beta\cdot\beta
=-\frac{\sin^2\theta}{\sin^2c}\,\exp\left(-2\int\lambda\cos\theta\,du\right),
\end{equation*}
Remark that the right hand sides are again constant and that $b$ is
a timelike vector, whereas $a$ and $\beta$ are spacelike. A suitable
isometry of $\H^2(-1)\times\R$, followed by the reparametrization
$v\sin c\mapsto v$, transforms the immersion given by (\ref{3.3.34})
and (\ref{3.3.23}) into \e{3.3.31}.

Finally, we solve the equations $\lambda^2-\sin^2\theta=-\sin^2c$,
$\theta'=\lambda\sin\theta$ explicitly. Putting $\theta=\arctan(f)$,
we obtain
$$\left(\frac{\theta'}{\sin\theta}\right)^2-\sin^2\theta=-\sin^2c \Leftrightarrow \frac{(f')^2}{f^2(f^2\cos^2c-\sin^2c)}=1.$$
We see that $f^2\cos^2c-\sin^2c>0$, and by integration, we obtain
$$\arctan\left(\frac{\sin
c}{\sqrt{f^2\cos^2c-\sin^2c}}\right)=\pm u\sin c+d \Leftrightarrow
f=\frac{\tan c}{\sin(\pm u\sin c+d)},$$ for some $d\in\R$. After a
translation in the $u$-coordinate, we obtain \e{3.3.31}. \hfill
$\square$\medskip

\begin{remark}
We can write the immersions of the first and the last case of
Theorem \ref{T:3.3.3} more explicitly. After the substitution
$e^{\pm cu}\mapsto u$, the immersion \e{3.3.27} becomes
\begin{equation*}
F(u,v)=\left(\frac{1+c^2+u^2}{p(u)q(u)},\ \frac{2u\cos
v}{p(u)q(u)},\ \frac{2u\sin v}{p(u)q(u)},\
\frac{1}{4c^2}\arctan\left(\frac{u^2-1+c^2}{2c}\right)\right),
\end{equation*}
with $p(u)=\sqrt{(u-1)^2+c^2}$ and $q(u)=\sqrt{(u+1)^2+c^2}$.

The immersion \e{3.3.31} is, after the substitution $\pm u\sin
c\mapsto u$, given by
\begin{equation*}\label{}
F(u,v)=\left(\frac{\cosh v}{p(u)},\ \frac{\sinh v}{p(u)},\
\frac{-\cos c\,\cos u}{p(u)},\ \arctan\left(\frac{\tan u}{\sin
c}\right)\right),
\end{equation*}
with $p(u)=\sqrt{1-\cos^2c\,\cos^2u}$.
\end{remark}

\begin{remark} Totally umbilical surfaces in BCV spaces and in the
Lie group $\Sol$ were independently studied in \cite{17aa}, from a
global viewpoint.
\end{remark}

\section{Higher order parallel surfaces}

The following example shows that every Hopf-cylinder in a BCV space
is flat.\medskip

\begin{example} Consider a Hopf-cylinder in $\M^3(\kappa,\tau)$. Let
$\{E_1=ae_1+be_2, E_2=e_3\}$, with $a^2+b^2=1$, be an orthonormal
frame field along the surface, then $N=E_1\times E_2=be_1-ae_2$ is a
unit normal. Using the equations in (\ref{Levi Civita}), one
computes
\begin{eqnarray*}
\nt_{E_1}N &=& \left(aE_1[b]-bE_1[a]+\frac{\kappa}{2}(ay-bx)\right)E_1-\tau E_2,\\
\nt_{E_2}N &=& (aE_2[b]-bE_2[a]-\tau)E_1.
\end{eqnarray*}
This means that the shape operator with respect to the basis
$\{E_1,E_2\}$ takes the form
\begin{eqnarray*}
S &=& \left(\begin{array}{cc} -aE_1[b]+bE_1[a]-\frac{\kappa}{2}(ay-bx) & -aE_2[b]+bE_2[a]+\tau \\ \tau & 0\end{array}\right)\\
  &=& \left(\begin{array}{cc} -aE_1[b]+bE_1[a]-\frac{\kappa}{2}(ay-bx) & \tau \\ \tau &
  0\end{array}\right),
\end{eqnarray*}
the last equation due to the symmetry. Remark that from this
symmetry we have $aE_2[b]=bE_2[a]$, which, together with $a^2+b^2=1$
implies that $a$ and $b$ are constant along the fibres of the
Hopf-fibration. From Gauss' equation (\ref{Gaussvgl2}), we have
$$K=\det S + \tau^2 + (\kappa-4\tau^2)\cos^2\theta = -\tau^2 +
\tau^2 + (\kappa-4\tau^2)\cos^2\frac{\pi}{2} = 0.$$
\end{example}

Now consider an arbitrary flat surface $M^2$ in $\M^3(\kappa,\tau)$.
Every $p\in M^2$ has an open neighbourhood $U$, which is isometric
to an open part of $\mathbb{E}^2$. Denote by $(u,v)$ the Euclidean
coordinates on $U$. Suppose $T=T_1\frac{\partial}{\partial
u}+T_2\frac{\partial}{\partial v}$ and $S=\left(S_{ij}\right)_{1\leq
i,j\leq 2}$ with respect to the orthonormal basis
$\left\{\frac{\partial}{\partial u}, \frac{\partial}{\partial
v}\right\}$. We consider $S_{11}$, $S_{12}$, $S_{22}$, $\cos\theta$,
$T_1$ and $T_2$ as functions of the Euclidean coordinates $(u,v)$ on
$U$.
\begin{lemma}\label{lemma3}
The functions $S_{11}$, $S_{12}$, $S_{22}$, $\cos\theta$, $T_1$ and
$T_2$ satisfy the following system of equations:
\begin{eqnarray}
T_1^2+T_2^2+\cos^2\theta=1;\label{1}\\
S_{11}S_{22}-S_{12}^2+\tau^2+(\kappa-4\tau^2)\cos^2\theta=0;\label{2}\\
\frac{\partial S_{12}}{\partial u}-\frac{\partial S_{11}}{\partial
v}=(\kappa-4\tau^2)T_2\cos\theta,\label{3}\\
\quad \frac{\partial S_{22}} {\partial u}-\frac{\partial S_{12}}{\partial v}=-(\kappa-4\tau^2)T_1\cos\theta;\nonumber\\
\frac{\partial T_1}{\partial u}=S_{11}\cos\theta, \quad
\frac{\partial T_1}{\partial v}=(S_{12}+\tau)\cos\theta,\label{4}\\
\frac{\partial T_2}{\partial u}=(S_{12}-\tau)\cos\theta,
\quad \frac{\partial T_2}{\partial v}=S_{22}\cos\theta;\nonumber\\
\frac{\partial\cos\theta}{\partial u}=-S_{11}T_1-S_{12}T_2+\tau T_2,\label{5}\\
\frac{\partial\cos\theta}{\partial v}=-S_{12}T_1-S_{22}T_2-\tau
T_1.\nonumber
\end{eqnarray}
\end{lemma}

\noindent\emph{Proof.} Equation (\ref{1}) follows immediately from
the definitions of $T$ and $\theta$. Equation (\ref{2}) expresses
Gauss' equation (\ref{Gaussvgl2}), while the equations (\ref{3})
express the equation of Codazzi (\ref{Codazzivgl}). The equations in
(\ref{4}) and (\ref{5}) are nothing but the structure equations
(\ref{structuurvgl1}) and (\ref{structuurvgl2}). \hfill $\square$
\medskip

The following result is the last step to obtain a full
classification of higher order parallel surfaces in BCV spaces.

\begin{proposition}\label{prop2}
A $k$-parallel, flat surface $M^2$ in a BCV space
$\M^3(\kappa,\tau)$, with $\kappa\neq4\tau^2$, is an open part of a
Hopf-cylinder over a curve in $\M^2(\kappa)$, whose curvature is a
polynomial function of degree at most $k-1$ of the arc length.
\end{proposition}

\noindent\emph{Proof.} Since $M^2$ is $k$-parallel and flat, the
functions $S_{11}$, $S_{12}$ and $S_{22}$ have to be polynomials of
degree at most $k-1$ in $u$ and $v$. First one can show that the
equations in lemma $\ref{lemma3}$ then imply that $\theta$ has to be
a constant. This proof is very similar to the proof of the Main
Theorem in \cite{11} and we will therefore omit it.

Now it follows from (\ref{3}) that the functions $T_1$ and $T_2$ are
polynomial functions in $u$ and $v$. Since $T_1$ and $T_2$ satisfy
$T_1^2+T_2^2=1-\cos^2\theta$ and $\theta$ is a constant, they have
to be constant. Then the equations in (\ref{4}) imply that either
$\cos\theta=0$ or $\tau=0$ and $S=0$. Totally geodesic surfaces in
BCV-spaces with $\tau=0$ are classified in proposition \ref{prop1}
and it is clear that the only flat ones are Hopf-cylinders. Hence we
may conclude that $M^2$ is an open part of a Hopf-cylinder.

To finish, we prove the assertion about the curvature of the base
curve. Taking $E_1$ and $E_2$ as in example 1, one can verify that
$\nabla_{E_i}E_j=0$ and hence we can take Euclidean coordinates
$(u,v)$ such that $E_1=\frac{\partial}{\partial u}$ and
$E_2=\frac{\partial}{\partial v}$. As we remarked before, $a$ and
$b$ will only depend on $u$ and we write $a'$ and $b'$ for the
derivatives with respect to $u$. The base curve
$\gamma(u)=(x(u),y(u))$ satisfies
$\gamma'=\pi_{\ast}E_1=(1+\frac{\kappa}{4}(x^2+y^2))(a,b)$, such
that $u$ is an arc length parameter. We compute
$$\kappa_{\gamma}=(1+\frac{\kappa}{4}(x^2+y^2))\frac{x'y''-x''y'}{((x')^2+(y')^2)^{\frac{3}{2}}}
+\frac{\kappa}{2}\frac{x'y-xy'}{((x')^2+(y')^2)^{\frac{1}{2}}}=ab'-a'b+\frac{\kappa}{2}(ay-bx)=-S_{11}.$$
Looking at the expression for $S$, we see that the surface is
$k$-parallel if and only if $S_{11}$ is a polynomial of degree at
most $k-1$ in $u$ and $v$. This is equivalent to $\kappa_{\gamma}$
being a polynomial of degree at most $k-1$ in $u$. \hfill $\square$
\medskip

From Lemma \ref{lemma1}, Proposition \ref{prop1} and Proposition
\ref{prop2} we obtain a full classification of higher order parallel
surfaces in 3-dimensional homogeneous spaces with 4-dimensional
isometry group:

\begin{theorem}\label{theo7} A $k$-parallel surface in a BCV space
$\widetilde{M}^3(\kappa,\tau)$, with $\kappa\neq4\tau^2$, is one of
the following:
\begin{itemize}
\item[(i)] an open part of a Hopf-cylinder on a curve whose geodesic curvature is a polynomial function of
      degree at most $k-1$ of the arc length;
\item[(ii)] an open part of a totally geodesic leaf of the Hopf-fibration;
\end{itemize}
the latter case only occuring when $\tau=0$.
\end{theorem}

\end{document}